\documentclass[a4paper,12pt]{elsarticle}

\usepackage{bm}

\usepackage{amssymb,amsxtra,graphicx} %amssymb,amsxtra - needed %amsmath, %amsthm, amsfonts,

\def\x#1{} %to mark changes
\def\xz{\hspace{-.17em}}
\def\xy{\hspace{.07em}}
\def\diag{\operatorname{diag}}
\def\Diag{\operatorname{Diag}}
\def\max{\operatorname{max}}
\def\R{{\mathbb R}}
\def\S{\Sigma}
\def\s{\kappa}
\def\bet{t}
\def\Comm{\operatorname{comm}}
\def\Katz{\operatorname{Katz}}
\def\df{\operatorname{df}}
\def\Heat{\operatorname{heat}}
\def\nHeat{\operatorname{n-heat}}
\def\reg{\operatorname{regL}}
\def\PPR{\operatorname{PPR}}
\def\ModifPPR{\operatorname{modifPPR}}
\def\absorp{\operatorname{absorp}}
\def\PPR{\operatorname{PPR}}
\def\heatPPR{\operatorname{heatPPR}}
\def\f{\varphi}
\def\cd{\cdot}
\def\De{D}
\def\D{{\cal D}}
\def\epr{\hfill$\square$\smallskip\par\medskip}\x{}

\newtheorem{thm}{Theorem}{\bfseries}{\itshape}
{\bfseries}{\itshape}
\newtheorem{cor}{Corollary}{\bfseries}{\itshape}
\newtheorem{pro}{Proposition}{\bfseries}{\itshape} %\x{}
\newtheorem{rem}{Remark}{\bfseries}{\upshape}

\begin{document}

\title{Similarities on Graphs: Kernels versus Proximity Measures{\large$^1$}\footnote{$\xz\!\!^1$This is the authors' version of a work that was submitted to Elsevier.\x{!}}}

\author[isa]{Konstantin Avrachenkov}\x{~}
\ead{k.avrachenkov@inria.fr}

%\author[ics,kot,mft]{Pavel Chebotarev\corref{cor1}}%\footnotemark[4]}
\author[ics]{Pavel Chebotarev\corref{cor1}}%\footnotemark[4]}
\ead{pavel4e@gmail.com}

\author[isa]{Dmytro Rubanov}\x{~}
\ead{dmytro.rubanov@inria.fr}

\address[isa]{Inria Sophia Antipolis, France}
\address[ics]{Trapeznikov\x{!} Institute of Control Sciences of the Russian Academy of Sciences, %Moscow, 
Russia}
%\address[kot]{Kotelnikov Institute of Radioengineering and Electronics of RAS, \x{!}Moscow, Russia}
%\address[mft]{Moscow Institute of Physics and Technology, \x{!}Moscow, Russia}
%RAS Institute of Control Sciences, Russia\\ {\tt pavel4e@gmail.com}}
 \cortext[cor1]{Corresponding author}

\begin{abstract}
%Kernels and, broadly speaking, similarity measures\x{} on graphs are extensively used in graph-based unsupervised and semi-supervised learning algorithms as well as in the\x{} link prediction problem.
We analytically study proximity and distance properties of various kernels and similarity measures on graphs.
This helps to understand the mathematical nature of such measures and can potentially be useful for recommending the adoption of specific similarity measures in data analysis.
\end{abstract}

\maketitle

\section{Introduction}

Until the 1960s, mathematicians studied only one distance for graph vertices, the shortest path distance~\cite{BH90}.
In 1967, Gerald Sharpe proposed the electric distance~\cite{S67}; then it was rediscovered several times.
For some collection of graph distances, we refer to~\cite[Chapter~15]{DD16}.

Distances\footnote{\x{!}In this paper, the term distance is used as a synonym of metric, i.e., every distance satisfies nonnegativity, symmetry, the identity of indiscernibles, and the triangle inequality.} are treated as dissimilarity measures.
In contrast, similarity measures are maximized when every distance is equal to zero, i.e., when two arguments of the function coincide.
At the same time, there is a close relation between distances and certain classes of similarity measures.

One of such classes consists of functions defined with the help of kernels on graphs, i.e., positive semidefinite matrices with indices corresponding to the nodes. Every kernel is the Gram matrix of some set of vectors in a Euclidean\x{!} space, and Schoenberg's theorem \cite{S35,S38} shows how it can be transformed into the matrix of Euclidean distances between these vectors.

Another class consists of proximity measures characterized by the triangle inequality for proximities. Every proximity measure $\s(x,y)$ generates a distance~\cite{CS98} by means of the $d(x,y)=\frac12(\s(x,x)+\s(y,y))-\s(x,y)$ transformation, and for some subclasses of proximity measures, the inverse transformation can also be specified.

Furthermore, it turns out that many similarity measures are transitional measures \cite{C11,C13}, in particular, they satisfy the inequality $s_{ij}\,s_{jk}\le s_{ik}\,s_{jj}$ and so they can be transformed into proximities by means of the logarithmic transformation.

Distances and similarity measures on graphs are widely used in data analysis, especially, in graph-based supervised, semi-supervised, and unsupervised machine learning, see, e.g., \cite{AMGS12,AGS13,ACM17,BL11,CSZ06,DGK04,FYPS06,FSS16,L98,LK07,MMRTS01,ZSH04} and references therein. In a number of studies including~\cite{IC17,KSS14,SFS16}, comparative ability of various measures to detect communities and predict links has been explored. However, in such studies, the authors do not focus on the mathematical properties of the measures under consideration.

The purpose of this paper is to start filling this gap. We consider a number of well-known and recently proposed similarity measures on graphs (including weighted graphs) defined in terms of one of the following basic matrices: the weighted adjacency matrix, the Laplacian matrix, and the (stochastic) Markov matrix. We explore their basic mathematical properties, in particular, we find out whether they belong to the classes of kernels or proximities and study the properties of distances related to them. This helps to reveal the nature of such measures and can be considered as a step towards a mathematical theory of similarity~/ dissimilarity measures on graphs.

\section{Definitions and preliminaries}

The \textit{weighted adjacency matrix\/} $W=(w_{ij})$ of a weighted undirected graph $G$ with
vertex set $V(G)=\{1,\ldots,n\}$ is the matrix with elements
$$
w_{ij}=\left\{\begin{array}{ll}
\mbox{weight of edge $(i,j)$}, & \mbox{if} \ i \sim j,\\
0, & \mbox{otherwise}.
\end{array}\right.
$$
In what follows, $G$ is connected.

The ordinary (or combinatorial) \textit{Laplacian matrix\/} $L$ of $G$ is defined as follows: $L=\De-W,$ where $\De=\Diag(W\!\cd\!\bm1)$ is the degree matrix of $G$, $\Diag(\bm x)$ is the diagonal matrix with vector $\bm x$ on the main diagonal, and $\bm1=(1,\ldots, 1)^T$. In most cases, the dimension of $\bm1$ is clear from the context.

Informally, given a weighted graph $G$, a \textit{similarity measure\/} on the set of its vertices $V(G)$ is a function $\s\!:V(G)\times V(G)\!\to\!\R$ that characterizes similarity (or affinity, or closeness) between the vertices of $G$ in a meaningful manner and thus is intuitively and practically adequate for empirical applications~\cite{AGS13,DFG01,FSS16,LK07}.

\textit{A kernel on graph\/} is a\x{!} graph similarity measure that has an inner product representation. All the\x{!} inner product matrices (also called Gram matrices) with real entries are symmetric positive semidefinite matrices. On the other hand, any semidefinite matrix has a representation as a Gram matrix with respect to the Euclidean inner product~\cite{HJ13}.

We note that following \cite{KL02,SK03} we prefer to write {\it kernel on graph} rather than {\it graph kernel},
as the notion of\x{} ``graph kernel'' refers to a kernel between graphs~\cite{Vetal10}.

A \textit{proximity measure\/} (or simply \textit{proximity}) \cite{CS98} on a %finite
set $A$ is a function $\s\!:A\!\times\! A\to\R$ that satisfies the
\textit{triangle inequality for proximities\/}, viz.:\\
For any $x,y,z\in A,$\; $\s(x,y)+\s(x,z)-\s(y,z)\le\s(x,x)$, and if $z=y$ and $y\ne x$, then the inequality is strict.

A proximity $\s$ is a \textit{$\S$-proximity\/} ($\S\in\R$) if it satisfies the \textit{normalization condition\/}$:$
$\sum_{y\in A}\s(x,y)=\S$ for any $x\in A.$

By setting $z = x$ in the triangle inequality for proximities and using the arbitrariness of $x$ and $y$ one verifies that any proximity satisfies \textit{symmetry}: $\s(x,y)=\s(y,x)$ for any $x,y\in A.$
Consequently, if $\s(x,y)$ is a proximity, then\, $-\s(x,y)$ is a protometric \cite{DC11,DD16}.

Furthermore, any $\S$-proximity has the \textit{egocentrism\/} property: $\s(x,x)>\s(x,y)$ for any distinct $x, y\in A$~\cite{CS98}.
If $\s(x,y)$ is represented by a matrix $K=(K_{xy})=(\s(x,y)),$ then egocentrism of $\s(x,y)$ amounts to the \textit{strict\x{"} entrywise diagonal dominance\/} of~$K.$

\x{\bm x}
If $\bm x_i$ and $\bm x_j$ are two points in the Euclidean space $\R^n,$ then $||\bm x_i-\bm x_j||_2^2$ is the squared distance between $\bm x_i$ and~$\bm x_j.$
Schoenberg's theorem  establishes a connection between positive semidefinite matrices (kernels) and matrices of Euclidean distances.
\begin{thm}[\cite{S35,S38}]
\label{t_S}
Let $K$ be an $n\times n$ symmetric matrix.
Define the matrix
\begin{equation}
\label{e_D2K}
\D = (d_{ij})\x{"} = \frac12\big(\diag(K)\cd\bm1^T + \bm1\cd\diag(K)^T\big)-K,
\end{equation}
where $\diag(K)$ is the vector consisting of the diagonal entries of~$K.$
Then there exists  a set of vectors $\bm x_1,\ldots,\bm x_n\in\R^n$ such that
$d_{ij} = ||\bm x_i-\bm x_j||_2^2$ $(i,j = 1,\ldots, n)$ if and only if $K$ is positive semidefinite.
\end{thm}

In the case described in Theorem~\ref{t_S}, $K$ is the Gram matrix of $\bm x_1,\ldots,\bm x_n$.
Given $K,$ these vectors can be obtained as the columns of the unique positive semidefinite real matrix $B$ such that $B^2=B^TB=K.$
$B$~has the expression $B=U\Lambda^{1/2}U^*,$ where $\Lambda=\Diag(\lambda_1,\ldots,\lambda_n)$, $\Lambda^{1/2}=\Diag(\lambda_1^{1/2},\ldots,\lambda_n^{1/2})$, and $A=U\Lambda U^*$ is the unitary decomposition of~$A$ \cite[Corollary~7.2.11]{HJ13}.

Connections between proximities and distances are established in~\cite{CS98}. %by Theorem~1
\begin{thm}%[\cite{CS98}]
\label{t_DS}
For any proximity $\s$ on a finite set $A,$ the function
\begin{equation}
\label{e_s2d}
d(x,y)=\frac12(\s(x,x)+\s(y,y))-\s(x,y),\quad x,y\in A
\end{equation}
is a distance function $A\!\times\! A\to\R.$\x{"}
%For any distance $d$ on a finite set $A$ with $|A|=n$, the function
%\begin{equation}
%\label{e_d2s}
%\s(x,y)=d(x,\cd)+d(y,\cd)-d(x,y)-d(\cd,\cd)+\frac{\S}n
%\end{equation}
\end{thm}

This theorem follows from the proof of Proposition~3 in~\cite{CS98}.

\begin{cor}
\label{c_sq}
Let $\D=(d_{xy})$ be obtained by \eqref{e_D2K} from a square matrix~$K.$ %of size~$n.$\x{}
If $\D$ has negative entries or $\sqrt{d_{xy}}+\sqrt{d_{yz}}<\sqrt{d_{xz}}$ for some $x,y,z\in\{1,\ldots,n\},$\x{} then the function $\s(x,y)=K_{xy},$\, $x,y\in\{1,\ldots,n\}$\x{} is not a proximity.
\end{cor}

{\bf Proof.} If $\sqrt{d_{xy}}+\sqrt{d_{yz}}<\sqrt{d_{xz}},$ then $d_{xy} + d_{yz} + 2\sqrt{d_{xy}d_{yz}}<d_{xz},$ i.e., the function $d(x,y)=d_{xy}$ violates the ordinary triangle inequality. Thus, it is not a distance, as well as in the case where $\D$ has negative entries. Hence, by Theorem~\ref{t_DS}, $\s$ is not a proximity.
\epr\x{}

The following theorem describes a one-to-one correspondence between distances and $\S$-proximities with a fixed $\S$ on the same finite set.

\begin{thm}[\cite{CS98}]
\label{t_DS1}
Let $\bm S$ and $\bm D$ be the set of\/ $\S$-proximities on $A$ $(|A|=n;$ $\S\in\R$ is fixed$)$ and the set of distances on $A,$ respectively.
Consider the mapping $\psi(\s)$ defined by \eqref{e_s2d}
%\begin{equation}
%\label{e_s2d}
%d(x,y)=\frac12(\s(x,x)+\s(y,y))-\s(x,y)
%\end{equation}
and the mapping $\f(d)$ defined by
\begin{equation}
\label{e_d2s}
\s(x,y)=d(x,\cd)+d(y,\cd)-d(x,y)-d(\cd,\cd)+\frac{\S}n,
\end{equation}
where $d(x,\cd)=\frac1n\sum_{y\in A}d(x,y)$ and $d(\cd,\cd)=\frac1{n^2}\sum_{y,z\in A}d(y,z).$
Then $\psi(\bm S)=\bm D,$\: $\f(\bm D)=\bm S,$ and $\f(\psi(\s)),\:\s\in\bm S$ and $\psi(\f(d)),\:d\in\bm D$ are identity transformations.
\end{thm}

\begin{rem}
\rm The $K\to \D$ transformation \eqref{e_D2K} is the matrix form of~\eqref{e_s2d}. %the mapping $\psi(\s)$ of Theorem~\ref{t_DS}.
The matrix form of \eqref{e_d2s} is
\begin{equation}
\label{e_HDH}
K=-H\D H+\S J,
\end{equation}
where $J=\tfrac1n\bm1\!\cdot\!\bm1^T$ and $H=I-J$ is the \textit{centering matrix}.
\end{rem}

\section{Kernel,\x{~} proximity, and distance properties}\label{s_pdp}

\subsection{Adjacency matrix based kernels}

Let us consider several kernels on graphs based on the weighted adjacency matrix $W$ of a graph.

\subsubsection{Katz kernel}

The {\it Katz kernel} \cite{K53} (also referred to as walk proximity \cite{CS98a} and von Neumann\footnote{M.~Saerens \cite{S17} has remarked that a more suitable name could be \textit{Neumann diffusion kernel}, referring to the \textit{Neumann series\/} $\sum_{k=0}^\infty T^k$ (where $T$ is an operator) named after Carl Gottfried Neumann, while a connection of that to John von Neumann is not obvious (the concept of von Neumann kernel in group theory is essentially different).} diffusion kernel \cite{KS-TC02,S-TC04}) is defined\footnote{In fact, L.~Katz considered $\sum_{k=1}^\infty (\alpha W)^k.$} as follows:
$$
K^{\Katz}(\alpha) = \sum_{k=0}^\infty (\alpha W)^k = [I-\alpha W]^{-1},
$$
with $0 < \alpha < (\rho(W))^{-1}$, where $\rho(W)$ is the spectral radius of $W$.

%Obviously, $K^{\Katz}$ is symmetric.
It is easy to see that $[I-\alpha W]$ is an M-matrix\footnote{For the properties of M-matrices, we refer to~\cite{KN12}.}, i.e., a matrix of the form $A = qI - B$, where $B=(b_{ij})$ with $b_{ij}\ge0$ for all $1\le i,j\le n,$ while $q$ exceeds the maximum of the moduli of the eigenvalues of~$B$ (in the present case, $q=1$). Thus, $[I-\alpha W]$ is a symmetric M-matrix, i.e., a Stieltjes matrix. Consequently, $[I-\alpha W]$ is positive definite and so is $K^{\Katz}(\alpha) = [I-\alpha W]^{-1}$. Thus, by Schoenberg's theorem, $K^{\Katz}$ can be transformed by \eqref{e_D2K} into a matrix of squared Euclidean distances.

Moreover, the Katz kernel has the following properties:\\
If $[I-\alpha W]$ is row diagonally dominant, i.e., $|1-\alpha w_{ii}|\ge\alpha\sum_{j\ne i}|w_{ij}|$ for all $i\in V(G)$ (by the finiteness of the underlying space, one can always choose $\alpha$ small enough such that
this inequality becomes valid) then
\begin{itemize}\x{"give strict proof}
\item $K^{\Katz}(\alpha)$ satisfies the triangle inequality for proximities (see Corollary~6.2.5 in~\cite{KN12}), therefore, transformation \eqref{e_s2d} provides a distance on $V(G)$;
\item $K^{\Katz}(\alpha)$ satisfies egocentrism (i.e., \textit{\x{"}strict entrywise\/} diagonal dominance;
see also Metzler's property in~\cite{KN12}).

\end{itemize}

%Since $W$ is symmetric for undirected graphs, we also have the symmetry property: $K^{\Katz}_{ij} = K^{\Katz}_{ji}.$
Thus, in the case of row diagonal dominance of $[I-\alpha W]$, the Katz kernel is a non-normalized %symmetric
proximity.

\subsubsection{Communicability kernel}

The {\it communicability kernel} \cite{FYPS06,EH07,EH08} is defined as follows:
$$
K^{\Comm}(t) = \exp(\bet W) = \sum_{k=0}^\infty \frac{\bet^k}{k!} W^k.
$$
(We shall use letter ``$t$'' whenever some notion of time can be attached to the kernel parameter;
otherwise, we shall keep using letter ``$\alpha$''.)
It is an instance of symmetric exponential diffusion kernels~\cite{KL02}.
Since $K^{\Comm}$ is positive semidefinite, by Schoenberg's theorem, it can be transformed by \eqref{e_D2K} into a matrix of squared Euclidean distances.
However, this does not imply that $K^{\Comm}$ is a proximity.

In fact, it is easy to verify that for the graph $G$ with weighted\x{!} adjacency matrix
\begin{equation}
\label{e_wpath}
W=\begin{pmatrix}
    0 & 2 & 0 & 0\\
    2 & 0 & 1 & 0\\
    0 & 1 & 0 & 2\\
    0 & 0 & 2 & 0
  \end{pmatrix},
\end{equation}
$K^{\Comm}(1)$ violates the triangle inequality for proximities on the triple of vertices $(1, \bm2, 3)$ (the ``$x$'' element of the inequality is given in bold). On the other hand, $K^{\Comm}(t)\to I$ as $t\to0,$ which implies that $K^{\Comm}(t)$ with a sufficiently small $t$ is a [non-normalized] proximity.

Note that the graph corresponding to \eqref{e_wpath} is a weighted path 1--2--3--4, and immediate intuition suggests the inequality $d(1,2)<d(1,3)<d(1,4)$ for a distance on its vertices. However, $K^{\Comm}(3)$ induces a Euclidean distance for which $d(1,3)>d(1,4).$ For $K^{\Comm}(4.5)$ we even have $d(1,2)>d(1,4).\,$
However, $K^{\Comm}(t)$ with a small enough positive $t$ satisfies the common intuition.

By the way, the Katz kernel behaves similarly: when $\alpha>0$ is sufficiently small, it holds that $d(1,2)<d(1,3)<d(1,4),$ but for $\alpha>0.375,$ we have $d(1,3)>d(1,4).$ Moreover, if $0.38795<\alpha<(\rho(W))^{-1},$\, then  $d(1,2)>d(1,4)$ is true.

\subsubsection{Double-factorial similarity}

The {\it double-factorial similarity} \cite{ES17} is defined as follows:
$$
K^{\df}(\bet) = \sum_{k=0}^\infty \frac{\bet^k}{k!!} W^k.
$$

As distinct from the communicability measure, $K^{\df}$ is not generally a kernel.
Say, for the graph with weighted adjacency matrix \eqref{e_wpath}, $K^{\df}(1)$ has two negative eigenvalues. Therefore $K^{\df}$ does not generally induce a set of points in $\R^n,$ nor does it induce a natural Euclidean distance on~$V(G).$

Furthermore, in this example, matrix $\D$ obtained from $K^{\df}(1)$ by \eqref{e_D2K} has negative entries. Therefore, by Corollary~\ref{c_sq}, the function $\s(x,y)=K_{xy}^{\df}(1),$\, $x,y\in V(G)$ is not a proximity.
%$K^{\df}(1)$ is not a proximity either. Say, for \eqref{e_wpath} it violates the triangle inequality for proximities on the triple of vertices $(1,2,3).$

However, as well as $K^{\Comm}(t),$\, $K^{\df}(t)\to I$ as $t\to0.$ Consequently, all eigenvalues of $K^{\df}(t)$ converge to~$1,$
and hence, $K^{\df}(t)$ with a sufficiently small positive $t$ satisfies the triangle inequality for proximities. Thus, $K^{\df}(t)$ with a small enough positive $t$ is a kernel and a [non-normalized] proximity.

\subsection{Laplacian based kernels}

\subsubsection{Heat kernel}

The {\it heat kernel} is a symmetric exponential diffusion kernel~\cite{KL02} defined as follows:
$$
K^{\Heat}(t) = \exp(-t L) = \sum_{k=0}^\infty \frac{(-t)^k}{k!} L^k,
$$
where $L$ is the ordinary Laplacian matrix of~$G.$
%Because $L$ is symmetric, so is the heat kernel.

$K^{\Heat}(t)$ is positive-definite for all values of $t$, and hence, it is a kernel.
Then, by Schoenberg's theorem, $K^{\Heat}$ induces a Euclidean distance on $V(G).$ For our example \eqref{e_wpath}, this distance
for all $t>0$ obeys the intuitive inequality $d(1,2)<d(1,3)<d(1,4).$
%However, this distance is rather specific. Say, for the graph with weighted adjacency matrix \eqref{e_wpath},

On the other hand, $K^{\Heat}$ is not generally a proximity. E.g., for the example \eqref{e_wpath}, $K^{\Heat}(t)$ violates the triangle inequality for proximities on the triple of vertices $(1,\bm2,3)$ whenever $t>0.431.$ As well as for the communicability kernel, $K^{\Heat}(t)$ with a small enough $t$ is a proximity. Moreover, it is a\x{!} 1-proximity, as \x{"}in the above series expansion, $L$ has row sums $0,$ while $L^0=I$ has row sums~$1.$ %and thus $K^{\Heat}(t)\bm1=\bm1.$
Thus, the 1-normalization condition is satisfied for any~$t>0$.\x{"}

\subsubsection{Normalized heat kernel}\x{~}

The {\it normalized heat kernel} is defined as follows:
$$
K^{\nHeat}(t) = \exp(-t {\cal L}) = \sum_{k=0}^\infty \frac{(-t)^k}{k!} {\cal L}^k,
$$
where ${\cal L}=\De^{-1/2}L\De^{-1/2}$ is the normalized Laplacian, $\De$ being the degree matrix of $G$~\cite{C97}.

For this kernel, the main conclusions are the same as for the standard heat kernel.
For the example \eqref{e_wpath}, $K^{\Heat}(t)$ violates the triangle inequality for proximities on the triple of vertices $(1,\bm2,3)$ when $t>1.497.$
It is curious to observe that the triangle inequality of the example \eqref{e_wpath} is violated \x{}starting with a larger value of $t$ in comparison with the case of the standard heat kernel.
An important distinction is that generally, ${\cal L}$ has nonzero row sums. As a result, $K^{\nHeat}$ does not satisfy the normalization condition, and even for small $t>0$,\x{"} $K^{\nHeat}$ is a non-normalized proximity.

\subsubsection{Regularized Laplacian kernel}\x{~}

The {\it regularized Laplacian kernel\/}, or {\it forest kernel\/} is defined \cite{CS95} as follows:
$$
K^{\reg}(\bet) = [I + \bet L]^{-1},
$$
where $t>0.$

As was shown in \cite{CS97,CS98a}, the regularized Laplacian kernel is a 1-proximity and a row stochastic matrix.
Since $[I + \bet L]$ is positive definite, so is $[I + \bet L]^{-1},$ and by Schoenberg's theorem, $K^{\reg}$ induces a Euclidean distance on~$V(G).$

For the example \eqref{e_wpath}, the induced distances corresponding to $K^{\reg}$ always satisfy $d(1,2)<d(1,3)<d(1,4).$ Regarding the other properties of $K^{\reg}$, we refer to~\cite{CS97,ACM17}.

It is the first encountered example of similarity measure that satisfies the both distance and proximity properties
for all values of the kernel parameter.

\subsubsection{Absorption kernel}\x{~}

The {\it absorption kernel} \cite{JT16} is defined as follows:
$$
K^{\absorp}(t) = [tA+L]^{-1},\quad t>0,\x{}
$$
where $A=\Diag(\bm a)$ and $\bm a=(a_1,\ldots,a_n)^T$\x{} is called the \textit{vector of absorption rates\/}\x{} and has positive components.
As $K^{\absorp}(t^{-1})=t(A+tL)^{-1},$ this kernel is actually a generalization of the previous one.

Since $[tA+L]$ is positive definite, Schoenberg's theorem attaches a matrix of squared Euclidean distances to $K^{\absorp}(t)$.

$[tA+L]$ is a row diagonally dominant Stieltjes matrix, hence, by Corollary~6.2.5 in~\cite{KN12} %Corollary 2.4 from \cite{CNX05}
we conclude that $K^{\absorp}$ satisfies the triangle inequality for proximities, i.e., $K^{\absorp}$ is a proximity (but not generally a $\S$-proximity).

\subsection{Markov matrix based kernels \x{}and measures}

\subsubsection{Personalized PageRank}

{\it Personalized PageRank} (PPR) similarity measure \cite{PBMW99} is defined as follows:
$$
K^{\PPR}(\alpha) = [I - \alpha P]^{-1},
$$
where $P=\De^{-1}W$ is a row stochastic (Markov) matrix, $\De$ is the degree matrix of $G,$ and $0<\alpha<1$, which corresponds to the
standard random walk on the graph.

In general, $K^{\PPR}(\alpha)$ is not symmetric, so it is not positive semidefinite, nor is it a proximity.

Moreover, the functions $d(x,y)$ obtained from $K^{\PPR}$ by transformation\footnote{If $K$ is symmetric, then \eqref{e_s2d2} coincides with \eqref{e_s2d}.}
\begin{equation}
\label{e_s2d2}
d(x,y)=\frac12(\s(x,x)+\s(y,y)-\s(x,y)-\s(y,x)) %,\quad x,y\in A
\end{equation}
need not generally be distances.
Say, for
\begin{equation}
\label{e_wpath5}
W=\begin{pmatrix}
    0 & 2 & 0 & 0 & 0\\
    2 & 0 & 1 & 0 & 0\\
    0 & 1 & 0 & 1 & 0\\
    0 & 0 & 1 & 0 & 2\\
    0 & 0 & 0 & 2 & 0
  \end{pmatrix}
\end{equation}
with $K^{\PPR}(\alpha),$ one has $d(1,3)+d(3,4)<d(1,4)$ whenever $\alpha>0.9515.$

$K^{\PPR}$ has only positive eigenvalues.
However, its symmetrized counterpart $\frac12(K^{\PPR}+(K^{\PPR})^T)$ may have a negative eigenvalue (say, with $\alpha\ge0.984$ for~\eqref{e_wpath} or with $\alpha\ge0.98$ for~\eqref{e_wpath5}). Thus, it need not be positive semidefinite and, consequently, by Theorem~\ref{t_S}, $\D$ obtained from it by \eqref{e_D2K} (or from $K^{\PPR}$ by~\eqref{e_s2d2}) is not generally a matrix of squared Euclidean distances.

\medskip
$K^{\PPR}$ satisfies the normalization condition.
For a\x{} small enough $\alpha$, it can be transformed (as well as $K^{\Comm}$ and $K^{\df}$) into a distance matrix using~\eqref{e_s2d2}.

\medskip
On the other hand, one can slightly modify Personalized PageRank so it becomes a proximity.
Rewrite $K^{\PPR}$\x{} as follows:
$$
[I - \alpha \De^{-1}W]^{-1} = [\De - \alpha W]^{-1} \De.\x{.}
$$
Then\x{",} consider

\subsubsection{Modified Personalized PageRank}

$$
K^{\ModifPPR}(\alpha) = [I - \alpha \De^{-1}W]^{-1} \De^{-1}
                      = [\De - \alpha W]^{-1},\quad 0<\alpha<1,
$$
which becomes a non-normalized proximity by Corollary~6.2.5 in~\cite{KN12}. In particular, the triangle inequality becomes
$$
\frac{K^{\PPR}_{ii}(\alpha)}{d_i} - \frac{K^{\PPR}_{ji}(\alpha)}{d_i} - \frac{K^{\PPR}_{ik}(\alpha)}{d_k} + \frac{K^{\PPR}_{jk}(\alpha)}{d_k} \ge 0,
$$
which looks like an interesting inequality for Personalized PageRank.
Due to symmetry, $K^{\ModifPPR}_{ij} = K^{\ModifPPR}_{ji}$, and we obtain an independent proof
of the following identity for Personalized PageRank \cite{AGS13}:
$$
\frac{K^{\PPR}_{ij}(\alpha)}{d_j} = \frac{K^{\PPR}_{ji}(\alpha)}{d_i}.
$$

Note that replacing the Laplacian matrix $L=\De-W$ with $\De-\alpha W$ is a kind of alternative regularization of~$L.$
Being diagonally dominant,
\begin{equation}\x{}
\label{e_bard}
\De-\alpha W=\bar dI-(\bar dI-\De+\alpha W)\,  %d_{\max} to \bar d
\end{equation}
(where $\bar d$\x{} is the maximum degree of the vertices of~$G$) is a Stieltjes matrix. Consequently, $\De-\alpha W$ is positive definite and so is $K^{\ModifPPR}(\alpha) = [\De-\alpha W]^{-1}$. Thus, by Schoenberg's theorem, $K^{\ModifPPR}$ can be transformed by \eqref{e_D2K} into a matrix of squared Euclidean distances.

\medskip
We note that Personalized PageRank can be generalized by using \x{\-}non-homo\-geneous restart \cite{AHS14},
which will lead to the discrete-time analog of the absorption kernel. However, curiously enough,
the discrete-time version has a smaller number of proximity--distance\x{!} properties than the continuous-time
version.

\subsubsection{PageRank heat similarity measure}

{\it PageRank heat  similarity measure} \cite{C07} is defined as follows:
$$
K^{\heatPPR}(t) = \exp(-t(I-P)).
$$
Basically, the properties of this measure are similar to those of the standard Personalized PageRank.
Say, for the example \eqref{e_wpath5} with $K^{\heatPPR},$ one has $d(1,2)+d(2,3)<d(1,3)$ whenever $t>1.45.$
%Since $K^{\heatPPR}$ is positive semidefinite, we conclude that $K^{\heatPPR}$ is a 1-proximity.

\subsection{Logarithmic similarity measures and transitional properties}\x{~}

Given a strictly positive similarity measure $s(x,y),$ the function $\s(x,y)=\ln s(x,y)$ is the corresponding \textit{logarithmic similarity\/}. %~\cite{C11}.

%a general framework was presented for constructing cutpoint additive distances. Namely,
Using Theorem~\ref{t_DS} it can be verified \cite{C11} that whenever $S=(s_{ij})=(s(i,j))$ produces a strictly positive \emph{transitional measure\/} on $G$ (i.e., $s_{ij}\,s_{jk}\le s_{ik}\,s_{jj}$ for all vertices $i$, $j$, and $k,$ while $s_{ij}\,s_{jk}=s_{ik}\,s_{jj}$ if and only if every path from $i$ to $k$ visits~$j$), we have\x{!} that the logarithmic similarity $\s(x,y)=\ln s(x,y)$ produces a \textit{cutpoint additive distance\/},
viz., a distance that satisfies $d(i,j)+d(j,k)=d(i,k)$ iff every path from $i$ to $k$ visits~$j$:\x{}
\begin{equation}
\label{e_cad}
d(i,j)=\tfrac12(\s(i,i)+\s(j,j)-\s(i,j)-\s(j,i))=\ln\sqrt{\frac{s(i,i)\xy s(j,j)}{s(i,j)\xy s(j,i)}}.
\end{equation}

In the case of digraphs, five transitional measures were indicated in \cite{C11}, namely,  \textit{connection reliability\/}, \textit{path accessibility\/} with a sufficiently small parameter, \textit{walk accessibility\/}, and two versions of \textit{forest accessibility\/}; the undirected counterparts of the two latter measures were studied in~\cite{C12}\x{} and~\cite{C11a}, respectively. %CheDeza12,CheBapatBalaji12

\x{}\begin{pro}
\label{p_TM}
$K^{\absorp},$ $K^{\PPR},$ and $K^{\ModifPPR}$ produce transitional measures.
\end{pro}

\x{}{\bf Proof.} For $K^{\absorp}(t) = [tA+L]^{-1},$ let $h=\max_i\{a_i t+d_i-w_{ii}\},$ where $d_i$ is the degree of vertex~$i.$ Then $K^{\absorp}(t)=[hI-(hI-tA-D+W)]^{-1}=[I-W']^{-1}h^{-1},$ where $W'=h^{-1}(hI-tA-D+W)$ is nonnegative with row sums less than~$1.$ Hence, $K^{\absorp}(t)$ is positively proportional to the matrix $[I-W']^{-1}$ of walk weights of the graph with weighted adjacency matrix~$W'.$

Similarly, by \eqref{e_bard},
$K^{\ModifPPR}(\alpha)= [\De-\alpha W]^{-1}=[I-W'']^{-1}\bar d^{-1}$, where $W''=\bar d^{-1}(\bar dI-\De+\alpha W)$ is nonnegative with row sums less than~$1.$
Consequently, $K^{\ModifPPR}(\alpha)$ is proportional to the matrix of walk weights of the graph whose weighted adjacency matrix is~$W''.$

Finally, $K^{\PPR}(\alpha)$ is the matrix of walk weights of the digraph with weighted adjacency matrix $\alpha P.$

Since by \cite[Theorem\:6]{C11}, any finite matrix of walk weights of a weighted digraph produces a transitional measure, so do $K^{\absorp},$ $K^{\PPR},$ and $K^{\ModifPPR}.$
\epr\x{}

\x{}Thus, as by Proposition~\ref{p_TM} and the results of \cite{C11}, $K^{\Katz},$ $K^{\reg},$ $K^{\absorp},$ $K^{\PPR},$ and $K^{\ModifPPR}$
produce transitional measures, we have that the corresponding \textit{logarithmic\/} dissimilarities \eqref{e_cad} are cutpoint additive distances.

\x{}Furthermore, if $S=(s_{ij})=(s(i,j))$ produces a strictly positive transitional measure on $G,$ then, obviously, $\s(x,y)=\ln s(x,y)$ satisfies $\s(y,x)+\s(x,z)-\s(y,z)\le\s(x,x),$ which coincides\footnote{On various alternative versions of the triangle inequality, we refer to~\cite{DC11}.} with the triangle inequality for proximities whenever $s(x,y)$ is symmetric. Therefore, as $K^{\Katz},$ $K^{\reg},$ $K^{\absorp},$ and $K^{\ModifPPR}$ are symmetric, we obtain that the corresponding logarithmic similarities $\s(x,y)=\ln s(x,y)$ are proximities.

\x{}$K^{\PPR}$ is not generally symmetric, however, it can be observed that $\widetilde K^{\PPR}$ such that $\widetilde K^{\PPR}_{ij}=\sqrt{K^{\PPR}_{ij}K^{\PPR}_{ji}}$ is symmetric and produces the same \textit{logarithmic\/} distance \eqref{e_cad} as $K^{\PPR}$.
Hence, the logarithmic similarity $\s(x,y)=\ln\widetilde K^{\PPR}_{xy}$ is a proximity.

\x{}%Consequently, it holds that the logarithmic similarities corresponding to \textit{regularized Laplacian kernel\/} and \textit{Katz kernel\/} are proximities. It is straightforward to show that this is also the case for the logarithmic measures corresponding to \textit{absorption kernel\/}, \textit{modified personalized PageRank\/}, and (after symmetrization) \textit{personalized PageRank\/}. %(and not the case for the other similarity measures under consideration).

\x{}At the same time, the above logarithmic similarities are not kernels, as the corresponding matrices have negative eigenvalues.

This implies that being a proximity is not a stronger property than being a kernel. By Corollary~\ref{c_sq}, the square root of the distance induced by a proximity is also a distance. However, this square rooted distance need not generally be Euclidean, thus, Theorem~\ref{t_S} is not sufficient to conclude that the initial proximity is a kernel.

It can be verified that all logarithmic measures corresponding to the similarity measures under study preserve the natural order of distances $d(1, 2) < d(1, 3) < d(1, 4)$ for the example~\eqref{e_wpath}.

\section{Conclusion}
We have considered similarity measures on graphs based upon three fundamental graph matrices:
the adjacency matrix, the Laplacian matrix, and the (stochastic) Markov matrix.
For each measure, we examine if it is a kernel, if it is a proximity measure, and if it is a transitional measure: these classes are not nested.
%The classification appears to be not straightforward and one cannot construct a simple hierarchical classification.
Regularized Laplacian turns out to be a similarity measure satisfying most of the properties, whereas the logarithmic similarity transformation appears to be useful as a tool for obtaining cutpoint additive distances.
We are currently working on an understanding of what consequences the established properties have for machine learning algorithms.

\section*{Acknowledgements}
The work of KA and DR was supported by the joint Bell Labs Inria ADR ``Network Science'' and
the work of PC was supported by the Russian Science Foundation (project no.16-11-00063 granted to IRE RAS).

%\section*{References}\x{!}

\end{document}